\documentclass{amsart}
\usepackage{amssymb}
\usepackage[all]{xy}
\usepackage{ifthen}

\newtheorem{thm}{Theorem}
\newtheorem{lemma}[thm]{Lemma}
\newtheorem{cor}[thm]{Corollary}
\newtheorem{propn}[thm]{Proposition}

\theoremstyle{remark}

\newtheorem{remark}[thm]{Remark}
\newtheorem{hypothesis}[thm]{\normalfont \bfseries Hypothesis}

\newtheorem{discussion}[thm]{Discussion}
\newtheorem*{conjHHSu}{\normalfont \bfseries Conjecture~(HHSu)}
\newtheorem*{conjHHSl}{\normalfont \bfseries Conjecture~(HHSl)}

\DeclareMathOperator{\height}{ht}
\DeclareMathOperator{\spec}{Spec}
\DeclareMathOperator{\depth}{depth}
\newcommand{\defeq}{:=}
\DeclareMathOperator{\tor}{Tor}
\DeclareMathOperator{\unm}{Unm}
\DeclareMathOperator{\ass}{Ass}
\newcommand{\naturals}{\mathbb{N}}
\DeclareMathOperator{\rhomo}{\widetilde H}
\DeclareMathOperator{\reg}{reg}
\DeclareMathOperator{\projdim}{pd}

\makeatletter

\@addtoreset{thm}{section}
\def\thethm{\thesection.\@arabic\c@thm}
\def\theenumi{\@alph\c@enumi}
\def\@oddfoot{\hfill\today}
\def\@evenfoot{\hfill\today}

\def\@lbibitem[#1]#2{\def\@biblab@l{#1}
      \ifthenelse{\equal{#2}{EGA}}{\def \@biblab@l{EGA}}{}
      \ifthenelse{\equal{#2}{M2}}{\def \@biblab@l{\texttt{M2}}}{}
      \item[\@biblabel{\@biblab@l}\hfill]\if@filesw
      {\let\protect\noexpand
       \immediate
       \write\@auxout{\string\bibcite{#2}{\@biblab@l}}}\fi\ignorespaces}

\makeatother

\title{Multiplicity Bounds for Quadratic Monomial Ideals}
\author{Manoj Kummini}
\address{University of Kansas\\Lawrence, KS 66045, USA.}
\email{kummini@math.ku.edu}

\begin{document}

\begin{abstract}
We prove the multiplicity bounds conjectured by Herzog-Huneke-Srinivasan
and Herzog-Srinivasan in the following cases: the strong conjecture for
edge ideals of bipartite graphs, and the weaker Taylor bound conjecture for
all quadratic monomial ideals. We attach a directed graph to a bipartite
graph with perfect matching, and describe operations on the directed graph
that would reduce the problem to a Cohen-Macaulay bipartite graph. We
determine when equality holds in the conjectured bound for edge ideals of
bipartite graphs, and verify that when equality holds, the resolution is
pure. We characterize bipartite graphs that have Cohen-Macaulay edge ideals
and quasi-pure resolutions.
\end{abstract}

\maketitle

\section{Introduction}
\label{sec:intro}

Let $V$ be a finite set, and let $R = \Bbbk[V]$ be a polynomial ring, over
a field $\Bbbk$, treating the elements of $V$ as indeterminates. We make
$R$ into a graded ring by setting $\deg x = 1$ for all $x \in V$. Let $f_1,
\cdots, f_m \in R$ be homogeneous polynomials, and let $I = (f_1, \cdots,
f_m)$. Set $c = \height I$. Let $e(R/I)$ denote the Hilbert-Samuel
multiplicity of $R/I$.

Let $\mathbb F_\bullet$ be a minimal graded free resolution of $R/I$ over
$R$.  Let $M_l \defeq M_l(I)$ be the largest twist with which $R$ appears
in $\mathbb F_l$, $1 \leq l \leq \projdim R/I$. Let $m_l \defeq m_l(I)$ be
the smallest twist with which $R$ appears in $\mathbb F_l$. These do not
depend on the choice of the resolution: since the $\tor_l^R(\Bbbk,R/I)$ are
graded, we can define the (graded) Betti numbers $\beta_{l,j}(R/I) =
\dim_\Bbbk \tor_l^R(\Bbbk,R/I)_j$. Then $m_l = \min\{j : \beta_{l,j}(R/I)
\neq 0\}$ and $M_l = \max\{j : \beta_{l,j}(R/I) \neq 0\}$. However the
$\tor_l^R(\Bbbk,R/I)$ are independent of the choice of the resolution of
$R/I$. Herzog-Huneke-Srinivasan~\cite{HeSrBnds98} conjectured that:
\makeatletter\protected@edef\@currentlabel{HHSu}\makeatother
\begin{conjHHSu}
\label{conj:HHSu}
For a homogeneous ideal $I$,
\[
e(R/I) \leq \frac{M_1M_2 \cdots M_c}{c!}.
\]
\end{conjHHSu}

This has subsequently been proved in various cases. A survey
appears in~\cite{FraSriMultConj07}. Some newer results
include~\cite{HeZhAsympt06, KubWelkbary06, MiroRoigDetIdealMult06,
NoSwCMconn06, PuthenMultConj06}.

Herzog-Huneke-Srinivasan further conjectured that:
\makeatletter\protected@edef\@currentlabel{HHSu}\makeatother
\begin{conjHHSl}
\label{conj:HHSl}
Assume that $R/I$ is Cohen-Macaulay. Then
\[
e(R/I) \geq \frac{m_1m_2 \cdots m_c}{c!}.
\]
\end{conjHHSl}

We say that $R/I$ \emph{has a pure resolution} if for each $l$, there is a
unique twist in the free module $\mathbb F_l$, or, equivalently, $M_l =
m_l$.  We say that $R/I$ \emph{has a quasi-pure resolution} if for each
$l$, $m_{l+1} \leq M_l$. Huneke-Miller~\cite{HuMiCMpure85} proved that if
$R/I$ is Cohen-Macaulay and has a pure resolution, then the above
conjectures hold, with equality.
Migliore-Nagel-R\"omer~\cite{MiNaRoExtns05} conjectured that:
\makeatletter\protected@edef\@currentlabel{HHSu}\makeatother
\begin{conjHHSl}
\label{conj:HHSe}
If equality holds in Conjecture~\eqref{conj:HHSu} or in
Conjecture~\eqref{conj:HHSl} then $R/I$ is Cohen-Macaulay with a pure
resolution.
\end{conjHHSl}

Herzog-Srinivasan~\cite{HeSrBnds98} proved that all the three conjectures
above are true when $R/I$ is Cohen-Macaulay and has a quasi-pure
resolution.

If additionally $f_1, \cdots, f_m$ are monomials, then $R/I$ has another
resolution $\mathbb T_\bullet$, called the Taylor resolution; see,
\textit{e.g.},~\cite[Ex. 17.11]{eiscommalg}. Let $T_l \defeq T_l(I)$ be the
largest twist with which $R$ appears in $\mathbb T_l$. Then $T_l = \max \{
\deg \mathrm{lcm}(f_{s_1}, \cdots, f_{s_l}) : 1 \leq s_1 < \cdots < s_l
\leq m\}$. Herzog-Srinivasan~\cite{HeSrMonBnds04} conjectured that:
\makeatletter\protected@edef\@currentlabel{HHSu}\makeatother
\begin{conjHHSl}
\label{conj:TB}
For a monomial ideal $I$,
\[
e(R/I) \leq \frac{T_1T_2 \cdots T_c}{c!}.
\]
\end{conjHHSl}

In general $T_l \geq M_l$ for all $1 \leq l \leq c$; hence
Conjecture~\eqref{conj:TB} is weaker than Conjecture~\eqref{conj:HHSu}. In
this paper we first prove Conjecture~\eqref{conj:TB} for all ideals generated by
quadratic monomials:

\def\thmquadTB{Let $I \subseteq R$ be generated by monomials of degree $2$.
Then
\[
e(R/I) \leq \frac{T_1T_2 \cdots T_c}{c!}.
\]}
{\begin{thm}
\label{thm:quadTB}
\thmquadTB
\end{thm}}

If $I$ is generated by square-free monomials of degree $2$, we define a
graph $G$ on $V$ by setting, for all $x, y \in V$, $xy$ to be an
edge of $G$ if and only if $xy$ is a minimal generator of $I$, and say that
$I$ is the \emph{edge ideal} of $G$. See~\cite[Chapter~6]{VillmonAlg01} for
a discussion on edge ideals. Note that such a graph $G$ is simple,
\textit{i.e.}, it has no loops and multiple edges between any pair of
vertices.

If $xy$ is an edge of $G$, then we say that $x$ and $y$ are
\emph{neighbours} of each other. An edge is \emph{incident} on its
vertices. We say that an edge $xy$ is \emph{isolated} if there are no other
edges incident on $x$ or on $y$. A vertex $x$ is a \emph{leaf vertex} if
there is a unique $y \in V$ such that $xy$ is an edge that is not isolated;
in this case, we call $y$ a \emph{stem vertex}, and refer to the edge $xy$
as a \emph{leaf}. The \emph{degree} of a vertex $x$, denoted $\deg_G x$,
is the number of edges incident on $x$. A \emph{tree} is a connected
acyclic graph, and a \emph{forest} is a graph in which each connected
component is a tree. A graph $G$ is \emph{bipartite}, if there is a
partition $V = V_1 \bigsqcup V_2$ and every edge of $G$ is of the form $xy$
where $x \in V_1$ and $y \in V_2$. (In this paper, $\bigsqcup$ denotes a
disjoint union.) Recall that a graph $G$ is bipartite if and only if it
does not contain odd cycles; in particular, forests are bipartite. A
\emph{path} is a tree in which every vertex has degree at most two. A
\emph{cycle} is a connected graph in which every vertex has degree exactly
two. We have that Conjecture~\eqref{conj:HHSu} holds for edge ideals of
bipartite graphs:
\def\thmquadHHSbipart{Let $I \subseteq R$ be the edge ideal of a bipartite
graph. Then
\[
e(R/I) \leq \frac{M_1M_2 \cdots M_c}{c!}.
\]}
\begin{thm}
\label{thm:quadHHSbipart}
\thmquadHHSbipart
\end{thm}

Let $I$ be the edge ideal of a graph $G$. We say that $G$ is
\emph{Cohen-Macaulay} if $R/I$ is Cohen-Macaulay.
Herzog-Hibi~\cite[Lemma~3.3 and Theorem~3.4]{HeHiCMbip05} (see
Theorem~\ref{thm:CMbipHeHi}, below) characterized Cohen-Macaulay bipartite
graphs.  \def\thmquadHHSebipart{Let $I$ be the edge ideal of a bipartite
graph. If equality holds in Conjecture~\eqref{conj:HHSu}, then $R/I$ is a
complete intersection, or is Cohen-Macaulay with $\reg R/I = 1$.  In either
of the cases, $R/I$ is Cohen-Macaulay and has a pure resolution.}
\begin{thm}
\label{thm:quadHHSe}
\thmquadHHSebipart
\end{thm}

Kubitzke-Welker~\cite{KubWelkbary06} proved that the
Conjecture~\eqref{conj:HHSu} holds for Stanley-Reisner ideals of barycentric
subdivisions of simplicial complexes; such ideals are square-free quadratic
monomial ideals, but most often they are not bipartite.
Novik-Swartz~\cite{NoSwCMconn06} establishes Conjecture~\eqref{conj:HHSu}
when $\dim R/I = n - c$ is small and studies the behaviour of $m_l$ in the
context of Conjecture~\eqref{conj:HHSl}.

The paper is organized as follows. In Section~\ref{sec:prelims} we make
some reductions. In Section~\ref{sec:taylor} we discuss some properties of
Taylor resolutions, and prove Theorem~\ref{thm:quadTB}.
Section~\ref{sec:rednToCM} is devoted to reducing the proof of
Conjecture~\eqref{conj:HHSu} to the Cohen-Macaulay case. In
Section~\ref{sec:cm-bipartite-resln}, we prove
Theorem~\ref{thm:quadHHSbipart}.  Finally, in
Section~\ref{sec:quasiPureRes}, following a proof of
Theorem~\ref{thm:quadHHSe}, we determine when a bipartite graph is
Cohen-Macaulay and its edge ideal has a quasi-pure resolution.

\section{Preliminary Remarks}
\label{sec:prelims}

Using polarization we first reduce to the case of $I$ generated by
square-free monomials; see~\cite[Exercise~3.15]{MiStCCA05} for details on
polarization. Herzog-Srinivasan~\cite{HeSrMonBnds04} showed that we can do
this in the context of Conjecture~\eqref{conj:TB}. In order to see that it
works for Conjecture~\eqref{conj:HHSu}, suppose that $\tilde I \subseteq
\tilde R$ is the polarization of $I$, in a larger ring $\tilde R$. Moving
from $\tilde I$ to $I$ preserves numerical data of a free resolution; in
particular, $\beta_{l,j}(\tilde R/\tilde I) = \beta_{l,j}(R/I)$, for all
$l,j$. On the other hand, the graded Betti numbers determine the Hilbert
function of $R/I$ and hence $e(R/I)$. Thus for the rest of the paper, we
assume that $I$ is square-free, \textit{i.e.}, $f_1, \cdots, f_m$ are
square-free monomials. For most part, this section does not assume that the
$f_l$ are quadratic.

If $\mathfrak p \subseteq R$ is a prime ideal such that $\height \mathfrak
p = c = \height I$ and $I \subseteq \mathfrak p$, then we say that
$\mathfrak p$ is an \emph{unmixed} prime of $R/I$. Denote the set of
unmixed primes of $R/I$ by $\unm R/I $. Unmixed primes are necessarily
minimal over $I$, so $\unm R/I  \subseteq \ass R/I$; we say that $I$ is
\emph{unmixed} if $\unm R/I  = \ass R/I$. Since $I$ is square-free,
$e(R/I)$ is the number of unmixed primes of $R/I$.  We add,
parenthetically, that for the edge ideal of a graph, unmixed primes
correspond to minimal vertex covers of size
$c$~\cite[Proposition~6.1.16]{VillmonAlg01}. 

We will need the following lemma several times in this paper:
\begin{lemma}[\protect{\cite[Proof of Theorem~$2.1$,
step~(i)]{HeZhAsympt06}}] 
\label{thm:MlincrHeZh}
Let $I \subseteq R$ be a homogeneous ideal. Then, for $2 \leq l \leq c$,
$M_l(I) > M_{l-1}(I)$.
\qed
\end{lemma}

We briefly discuss multigraded resolutions and Hochster's
formula. By $\naturals^V$, we denote the set of functions from $V$ to the
set of natural numbers, $\naturals$, and by $\mathbf e_x$, the function
that sends $x \mapsto 1$ and $y \mapsto 0$ for all $y \neq x$, for all $x,
y \in V$.  We treat $R$ as $\naturals^V$-graded, by setting, for all $x \in
V$, $\deg x = \mathbf e_x$. When an $R$-module is $\naturals^V$-graded, we
will say that it is \emph{multigraded}. Since $I$ is a monomial ideal,
$R/I$ is multigraded, and so are the $\tor_l^R(\Bbbk,R/I)$. A
\emph{multidegree} $\sigma$ is an element of $\naturals^V$. We can then
define \emph{multigraded Betti numbers} $\beta_{l,\sigma}(R/I) \defeq
\dim_\Bbbk \tor_l^R(\Bbbk,R/I)_\sigma$. If $\sigma \in \naturals^V$, we
also represent the multidegree $\sigma$ as the monomial $\prod_{x \in
V}x^{\sigma(x)}$.

(We will use the same notation $\beta_{l,\cdot}$ for Betti numbers in both
the standard grading and multigrading. Notice that $\beta_{l,j}(\cdot) =
\sum\beta_{l, \sigma}(\cdot)$, where the sum is taken over the set
of $\sigma$ such that $\sum_{x \in V} \sigma(x) = j$.)

Let $\Delta$ be the Stanley-Reisner complex of $I$. The correspondence
between non-faces of $\Delta$ and monomials in $I$ can also be expressed as
follows: for any monomial prime ideal $\mathfrak p \in \spec R$, $I
\subseteq \mathfrak p$ if and only if $\mathfrak p = (\bar F)R$, the ideal
generated by $\bar F \defeq V \setminus F$, for some $F \in
\Delta$~\cite[Theorem~1.7]{MiStCCA05}. Thus, minimal prime ideals of $R/I$
correspond to complements of maximal faces of $\Delta$. If $I$ is
the edge ideal of a graph $G$, then $\Delta$ is the clique complex of the
complement graph $\bar G$~\cite[Section~6.2]{VillmonAlg01}. In this case,
we say that $\Delta$ is the \emph{coclique complex} of $G$. 

Proposition~\ref{thm:hochsterFormula} below will show that most
multidegrees of interest in this paper are square-free, \textit{i.e.},
$\sigma(x) \in \{0,1\}$ for all $x \in V$. When a multidegree
$\sigma$ is square-free, we will also use subsets of $V$ to represent
multidegrees, by identifying $\sigma \in \naturals^V$ with $\prod_{x :
\sigma(x) = 1} x$ and with $\{x \in V :
\sigma(x) = 1\}$. For $\sigma \subseteq V$, $\Delta|_\sigma$ denotes the
simplicial complex obtained by taking all the faces of $\Delta$ whose
vertices belong to $\sigma$. Similarly, we define $G|_\sigma$ to be the the
subgraph \emph{induced} on $\sigma$, \textit{i.e.}, the subgraph of $G$
obtained by taking all the edges between pairs of vertices in $\sigma$.

\begin{propn}[Hochster's Formula \protect{\cite[Corollary~5.12]{MiStCCA05}}]
\label{thm:hochsterFormula}
Suppose that $I$ is a square-free monomial ideal. Let $\Delta$ be the
Stanley-Reisner complex of $I$. The multidegrees $\sigma$ with
$\beta_{l,\sigma}(R/I) \neq 0$ are square-free, and
\[
\beta_{l,\sigma} (R/I) = 
\rhomo_{\vert \sigma \vert - l - 1}(\Delta \vert_\sigma; \Bbbk).
\]
\end{propn}
Here $\vert \cdot \vert$ denotes the cardinality of a set; later, in
Section~\ref{sec:rednToCM}, we will also use it to denote the underlying
undirected graph of a directed graph. We now describe how the graded Betti
numbers change under restriction to a subset of the variables and under
taking colons.
\begin{lemma}
\label{thm:varSetRestrictColon}
Let $I \subseteq R = \Bbbk[V]$ be a square-free monomial ideal, $x \in
V$, $l, j \in \naturals$ and $\sigma \subseteq V$ with $|\sigma| = j$. Then
\begin{enumerate}

\item \label{thm:varSetRestrict} Let $W \subseteq V$ and $J = (I \cap
\Bbbk[W])R$. Then,
\[
\beta_{l,\sigma}(R/J) = 
\begin{cases}
0, & \sigma \nsubseteq W, \\
\beta_{l,\sigma}(R/I), & \sigma \subseteq W.
\end{cases}
\]
In particular, $\beta_{l,j}(R/J) \leq \beta_{l,j}(R/I)$.

\item \label{thm:varColon} If $\beta_{l,\sigma}(R/(I:x)) \neq 0$, then
$\beta_{l,\sigma}(R/I) \neq 0$ or $\beta_{l,\sigma\cup\{x\}}(R/I) \neq
0$.
\end{enumerate}
\end{lemma}

\begin{proof}
\eqref{thm:varSetRestrict}: The second assertion follows from the first,
which we now prove. Let $\tilde \Delta$ be the Stanley-Reisner complex of
$J$.  Since for all $x \in V \setminus W$, $x$ does not belong to any
minimal prime ideal of $R/J$, we see that every maximal face of $\tilde
\Delta$ is contains $V \setminus W$. Hence if $\sigma \not \subseteq W$,
then for all $x \in \sigma \setminus W$, $\tilde \Delta |_\sigma$ is a cone
with vertex $x$, which, being contractible, does not have any homology.
Applying Proposition~\ref{thm:hochsterFormula}, we see that
$\beta_{l,\sigma}(R/J) = 0$.

Now let $\sigma \subseteq W$ and $F \subseteq V$. Then $F \in \Delta
\vert_\sigma$ if and only if $I \subseteq (\bar F)R$ and $F \subseteq
\sigma$ if and only if $J \subseteq (\bar F)R$ and $F \subseteq \sigma$ if
and only if $F \in \tilde \Delta \vert_\sigma$. Apply
Proposition~\ref{thm:hochsterFormula} again to get \[
\beta_{l,\sigma} (R/J) = 
\rhomo_{\vert \sigma \vert - l - 1}(\tilde \Delta \vert_\sigma; \Bbbk)
= \rhomo_{\vert \sigma \vert - l - 1}(\Delta \vert_\sigma; \Bbbk)
= \beta_{l,\sigma} (R/I).
\]

\eqref{thm:varColon}: We take the multigraded exact sequence of
$R$-modules:
\begin{equation}
\label{seq:fundColSpecialMulti}
\xymatrix{%
0 \ar[r] & \frac{R}{(I:x)}(-x) \ar[r] & \frac{R}I \ar[r] &
\frac{R}{(I,x)} \ar[r] & 0.}
\end{equation}
The corresponding multigraded long exact sequence of $\tor$ is
\[
\xymatrix{%
\cdots \ar[r] & \tor_{l+1}(\Bbbk, \frac{R}{(I,x)}) \ar[r] &
\tor_l(\Bbbk, \frac{R}{(I:x)}(-x)) \ar[r] & 
\tor_l(\Bbbk, \frac{R}I) \ar[r]
& \cdots.}
\]

Let $W = V \setminus \{x\}$ and $J = (I \cap \Bbbk[W])R$. Since
$\beta_{l,\sigma}(R/(I:x)) \neq 0$ and $x$ does not divide any monomial
minimal generator of $(I:x)$, we have, by the same argument as in
\eqref{thm:varSetRestrict}, $\sigma \subseteq W$. Let $\tau = \sigma
\cup \{x\}$. First observe that 
\[
\tor_l\left(\Bbbk, \frac{R}{(I:x)}\right)_\sigma \simeq 
\tor_l\left(\Bbbk, \frac{R}{(I:x)}(-x)\right)_\tau.
\]
Let us assume that $\beta_{l,\tau}(R/I) = 0$, because, if
$\beta_{l,\tau}(R/I) \neq 0$, there is nothing to prove. Restricting the
above long exact sequence to multidegree $\tau$, we see that
$\tor_{l+1}(\Bbbk, \frac{R}{(I,x)}) \neq 0$. Now , since $(I,x) =
(J,x)$, we see further $\tor_{l+1}(\Bbbk, \frac{R}{(J,x)})_\tau
\neq 0$. 

Since $x$ is a non-zerodivisor on $R/J$, we have a multigraded short
exact sequence
\[
\xymatrix{%
0 \ar[r] & \frac{R}J(-x) \ar[r] & \frac{R}J \ar[r] &
\frac{R}{(J,x)} \ar[r] & 0,}
\]
which gives the following long exact sequence of
$\tor$:
\[
\xymatrix{%
\cdots \ar[r] & 
\tor_{l+1}(\Bbbk, \frac{R}J) \ar[r] &
\tor_{l+1}(\Bbbk, \frac{R}{(J,x)}) \ar[r] &
\tor_l(\Bbbk, \frac{R}J(-x)) \ar[r] & \cdots.}
\]

Since $x$ does not divide any minimal monomial generator of $J$,
$\beta_{l+1, \tau}(R/J) = 0$. Therefore $\tor_l(\Bbbk, \frac{R}J(-x))_\tau
\neq 0$, or, equivalently, $\tor_l(\Bbbk, \frac{R}J)_\sigma \neq 0$.
By~\eqref{thm:varSetRestrict} above, $\beta_{l,\sigma} (R/I) \neq 0$.
\end{proof}

\begin{remark}
\label{rmk:ColonHasSameDepth}
Lyubeznik showed that, with notation as above, $\depth R/(I:x) \geq \depth
R/I$~\cite[Lemma~1.1]{LyubNonCMMon88};
Lemma~\ref{thm:varSetRestrictColon}\eqref{thm:varColon} gives another
proof.
\end{remark}

\begin{propn}
\label{thm:MlTlcolons}
Let $1 \leq l \leq c$. Then for all $x \in V$, 
\begin{enumerate}
\item \label{thm:MlDecr}
$M_l((I,x)) \leq M_l(I)$ and  $M_l((I:x)) \leq M_l(I)$. 
\item \label{thm:TlDecr}
$T_l((I,x)) \leq T_l(I)$ and $T_l((I:x)) \leq T_l(I)$. 
\end{enumerate}
\end{propn}

\begin{proof} Let $W = V \setminus \{x\}$ and $J = (I \cap \Bbbk[W])R$.
Then $(I,x) = (J,x)$ and $x$ is a non-zerodivisor on $R/J$; hence $c-1 \leq
\height J \leq c$.

\eqref{thm:MlDecr}: Let $\mathbb G_\bullet$ be a minimal graded
free resolution of $\mathrm R/J$. Denote the (graded) Koszul complex on
$x$ by $\mathbb K_\bullet$. Then $\mathbb G_\bullet \otimes \mathbb
K_\bullet$ is a minimal graded free resolution of $R/((J,x))$; in
particular, $M_l((I,x)) = M_l((J,x)) = \max \{ M_l(J), M_{l-1}(J) +
1 \}$. Since $\height J \geq c-1$, we conclude using
Lemmas~\ref{thm:MlincrHeZh}
and~\ref{thm:varSetRestrictColon}\eqref{thm:varSetRestrict} that, for $1
\leq l \leq c-1$, $M_l((I,x)) = M_l(J) \leq M_l(I)$. If $M_c(J) >
M_{c-1}(J)$, then $M_c((I,x)) = M_c(J) \leq M_c(I)$; otherwise, $M_c((I,x))
= M_{c-1}(J) + 1 \leq M_{c-1}(I) + 1 \leq M_c(I)$.

Lemma~\ref{thm:varSetRestrictColon}\eqref{thm:varColon} implies that
$M_l((I:x)) \leq M_l(I)$.

\eqref{thm:TlDecr}: Recall that $I$ is generated by square-free monomials
$f_1, \cdots, f_m$. Let
\[(f_j:x)
\defeq 
\begin{cases}
\frac{f_j}{x}, & \ \text{if} \ x \ \text{divides} \ f_j, \\
f_j, & \text{otherwise}
\end{cases}
\]
Since $(I:x) = ((f_1:x), \cdots, (f_m:x))$ and $(I,x) = (J,x)$,
the conclusions follow easily from the definition of $T_l$.
\end{proof}

Consider $I \cap R_1$, the vector space generated by the linear forms in
$I$.  Suppose that $\dim_\Bbbk (I \cap R_1) > 0$; then, since $I$ is a
monomial ideal, there exists $x \in V$ such that $x \in I$.
Write $J = (I \cap \Bbbk[V \setminus \{x\}])R$. Then $\height J = c-1$ and
$I = (J, x)$. Note that $e(R/J) = e(R/I)$. From
Lemma~\ref{thm:varSetRestrictColon}\eqref{thm:varSetRestrict} we know that
$M_l(J) \leq M_l(I)$ for $1 \leq l \leq c-1$. From the definition of $T_l$,
we see that $T_l(J) \leq T_l(I)$ for $1 \leq l \leq c-1$. Therefore it is
enough to prove Conjectures~\eqref{conj:HHSl} and~\eqref{conj:TB} for $J$.
In other words, $I$ behaves like an ideal of height $c-1$. Hence, if $\dim_\Bbbk (I \cap R_1) = \delta$, we will say that $I$ is
\emph{essentially of height $c-\delta$}.

\begin{discussion}
To make further reduction, we use the
sequence~\eqref{seq:fundColSpecialMulti}. Let $x \in V$. If $\height (I:x)
> c$, then $e(R/I) = e(R/(I,x))$. In light of
Proposition~\ref{thm:MlTlcolons}, we
can replace $I$ by $I$ by $(I,x)$ which is essentially of height $\leq
c-1$, and prove Conjectures~\eqref{conj:HHSu} and~\eqref{conj:TB} by
induction on height. We can also look at $(I,x)$ as an ideal in $n-1$
variables. On the other hand, if $\height (I,x) > c$, then
$e(R/I) = e(R/(I:x))$; we then replace $I$ by $(I:x)$ which is an ideal in
$n-1$ variables. In this case, we can prove the conjectures using induction
on the number of variables. Therefore, we reduce to the case that $\height
(I:x) = c = \height (I,x)$. For later use, we record this below:
%
\begin{hypothesis}
\label{standing-hyp}
For all $x \in V$, $\height (I:x) = c =
\height (I,x)$; consequently, $e(R/I) = e(R/(I,x)) + e(R/(I:x))$.
\end{hypothesis}

The remark about $e(R/I)$ follows from~\eqref{seq:fundColSpecialMulti}.
This hypothesis is equivalent to the assumption that for all $x \in V$,
there exist $\mathfrak p, \mathfrak q \in \unm R/I$ such that $x \in
\mathfrak p \setminus \mathfrak q$. Moreover, while proving the
conjectures, we will assume, inductively, that
conjectures~\eqref{conj:HHSu} and~\eqref{conj:TB} hold for $(I:x)$ and
$(I,x)$. 

We now restrict our attention to (square-free) quadratic monomial ideals,
\textit{i.e.}, $f_1, \cdots, f_m$ are square-free monomials of degree $2$.
Since $I$ is square-free, $I$ is the edge ideal of a graph $G$. For any $x
\in V$, $(I,x)$ is essentially of height $\leq c-1$, and $(I:x)$ is
essentially of height $\leq c-\delta$, where $\delta$ is the degree of the
vertex $x$ in $G$. If $G$ is bipartite, then
Hypothesis~\ref{standing-hyp} means that $G$ is perfectly matched; see
Section~\ref{sec:rednToCM}. 
\end{discussion}

\begin{discussion}
\label{disc:defn-of-mu-r-c}
For $\rho, \gamma \in \naturals$, let
\[
\mu(\rho,\gamma) \defeq 
\begin{cases}
\frac{(2\rho+1) \cdots (\gamma+\rho)}{(\rho+1)\cdots \gamma}, 
& \rho<\gamma \\
1, & \text{otherwise}.
\end{cases}
\]
Note that $\mu(\gamma-1,\gamma) = (2\gamma-1)/\gamma < 2 = 2\mu(\gamma,
\gamma)$. For any $\rho <
\gamma-1$, 
\[
\frac{\mu(\rho+1, \gamma)}{\mu(\rho, \gamma)} = 
\frac{(\rho+1)(\gamma+\rho+1)}{(2\rho+1)(2\rho+2)} > \frac12.
\]
Combining these, we conclude that 
\begin{equation}
\label{eqn:ratios-of-mu-r-c}
2\mu(\rho+1,\gamma) > \mu(\rho,\gamma), \; \text{for all}\; \rho <
\gamma \in \naturals. 
\end{equation}
\end{discussion}

We conclude this section with a crucial lemma that captures the main
numerical argument in the proofs of Theorems~\ref{thm:quadTB}
and~\ref{thm:quadHHSbipart}.
\begin{lemma}
\label{thm:mu-r-c-c1-numer-crit}
Let $\rho, \gamma, \gamma_1 \in \naturals$ such that $2 \leq \rho < \gamma
\leq \rho\gamma_1$ and $\rho-1 \leq \gamma-\gamma_1$. Then 
\[
2^\rho\mu(\rho,\gamma-1)+2^{\rho-1}\mu(\rho-1,\gamma-\gamma_1) <
2^\rho\mu(\rho,\gamma).
\]
\end{lemma}

\begin{proof}
Since $\mu(\rho,\gamma) - \mu(\rho,\gamma-1) =
\rho\mu(\rho,\gamma-1)/\gamma$, we must show that
\begin{equation}
\label{eqn:muRhoMinusOneGamma}
\mu(\rho-1,\gamma-\gamma_1) < \frac{2\rho\mu(\rho,\gamma-1)}\gamma.
\end{equation}

We first reduce the problem to the case of $\gamma=\rho\gamma_1$ as
follows. If $\gamma < \rho\gamma_1$, and if we replace $\gamma$ by
$\gamma+1$, the left-hand-side and right-hand-side
of~\eqref{eqn:muRhoMinusOneGamma} change by factors of
\[
\frac{\gamma-\gamma_1+\rho}{\gamma-\gamma_1+1} \quad
\text{and} \quad
\frac{\gamma+\rho}{\gamma+1}
\]
respectively. Both these factors are greater than $1$, and the
left-hand-side increases by a larger factor than the right-hand-side.
Therefore, it is enough to prove the lemma when $\gamma=\rho\gamma_1$,
\textit{i.e.}, that 
\[
\mu(\rho-1,\rho\gamma_1 -\gamma_1) <
\frac{2\rho\mu(\rho,\rho\gamma_1-1)}{\rho\gamma_1}.
\]
The hypothesis gives that $\gamma_1 > 1$, so we need to show that
\[
\frac{(2\rho-1)(2\rho)\cdots(\rho\gamma_1-\gamma_1+\rho-1)}{\rho(\rho+1)\cdots
(\rho\gamma_1-\gamma_1)} < 
\frac{2\rho}{\rho\gamma_1}
\frac{(2\rho+1)(2\rho+2)\cdots(\rho\gamma_1+\rho-1)}{(\rho+1)(\rho+2)\cdots
(\rho\gamma_1-1)}.
\]
We can verify this by hand for $(\rho,\gamma_1) = (2,2), (2,3)$ and
$(3,2)$. For all other values of $\rho,\gamma_1$, $\rho+1 \leq
\rho\gamma_1-\gamma_1-1$ and we rewrite the above equation as
\[
\frac{(2\rho-1)}\rho \cdot 2\rho \cdot \frac{(2\rho+1)}{(\rho+1)} \cdots
\frac{(\rho\gamma_1-\gamma_1-1+\rho)}{(\rho\gamma_1-\gamma_1-1)}
\frac1{(\rho\gamma_1-\gamma_1)} < 
\frac{2\rho}{\rho\gamma_1} \frac{(2\rho+1)}{(\rho+1)} \cdots
\frac{(\rho\gamma_1-1+\rho)}{(\rho\gamma_1-1)}
\]
which is equivalent to the following sequence of equivalent statements:
\begin{align*}
\frac{(2\rho-1)}\rho \frac1{(\rho\gamma_1-\gamma_1)} & <
\frac{(\rho\gamma_1-\gamma_1+\rho)}{(\rho\gamma_1-\gamma_1)} \cdots
\frac{(\rho\gamma_1-1+\rho)}{(\rho\gamma_1-1)} \frac1{\rho\gamma_1} \\
\frac{(2\rho-1)}\rho & <
\frac{(\rho\gamma_1-\gamma_1+\rho)}{(\rho\gamma_1-\gamma_1+1)} \cdots
\frac{(\rho\gamma_1-1+\rho)}{\rho\gamma_1} \\
\left(1 + \frac{\rho-1}\rho\right) & <
\left(1 + \frac{\rho-1}{(\rho\gamma_1-\gamma_1+1)}\right) \cdots 
\left(1 + \frac{\rho-1}{\rho\gamma_1}\right).
\end{align*}
This is indeed true, since there are $\gamma_1$ terms on the
right-hand-side and each of them is at least as large as $\left(1 +
\frac{\rho-1}{\rho\gamma_1}\right)$. Recall that $\gamma_1>1$.
\end{proof}

\section{Proof of Theorem~\protect{\ref{thm:quadTB}}}
\label{sec:taylor}

We first make some observations on how $T_l$ changes with $l$. Let
$\rho(I)$ be the length of the longest $R$-regular sequence in $\{f_1,
\cdots, f_m\}$. 

\begin{lemma}
\label{thm:first-two-then-one}
Assume Hypothesis~\ref{standing-hyp}. Then, for all $1 \leq l \leq m$, if
$T_l < n$, then $T_l > T_{l-1}$.  For all $2 \leq l \leq m-1$, we have $T_l
- T_{l-1} \geq T_{l+1} - T_l$. Consequently,
\[
T_l = 
\begin{cases}
2l, & 1 \leq l \leq \rho(I) \\
\min\{\rho(I) + l, n\}& \rho(I) \leq l \leq c
\end{cases}
\]
Moreover, for all $x \in V$, $\rho((I:x)) < \rho(I)$. 
\end{lemma}

\begin{proof}
A consequence of Hypothesis~\ref{standing-hyp} is that for every $x \in V$ there is a
monomial minimal generator $f_j$ such that $x$ divides $f_j$, from which
the first assertion follows. To prove the second assertion, assume, by way
of contradiction, and by induction on $m$, that $m$ is the smallest integer
$m'$ such that there exists an ideal generated by $m'$ quadratic monomials
such that the conclusion does not hold. Write $\delta_l = T_l - T_{l-1}$;
it is clear that $0 \leq \delta_l \leq 2$. Pick $l$ smallest such that
$\delta_l < \delta_{l+1}$. If $\delta_l = 0$, then $T_{l+1} = T_l = n$.
Hence $\delta_l = 1$ and $\delta_{l+1} = 2$.

We now claim that $l = m-1$. For, assume, without loss of generality, that
$T_{l+1} = \deg \mathrm{lcm}(f_1, \cdots, f_{l+1})$. Let $J =(f_1, \cdots,
f_{l+1})$. Then $T_l(J) \leq T_l(I) = T_{l+1}(I) - 2$. If $m > l+1$, then,
by minimality of $m$, $\delta_2(J) = \cdots = \delta_{l+1}(J) = 2$. Hence
$f_1, \cdots, f_{l+1}$ is a regular sequence, and, therefore, $T_j(I) = 2j,
\forall j \leq l+1$ and $\delta_2(I) = \cdots = \delta_{l+1}(I) = 2$
contradicting the choice of $l$. Therefore $l = m-1$.

Assume that $T_{m-1}(I) = \deg \mathrm{lcm} (f_1, \cdots, f_{m-1})$.  Let
$J = (f_1, \cdots, f_{m-1})$. If $T_{m-2}(J) < T_{m-2}(I)$, then
$\delta_{m-1}(J) = 2$, and hence $J$ is generated by a regular sequence of
$m-1$ quadratic monomials. Therefore $T_j(J) = 2j, 1 \leq j \leq m-1$. Since $T_j(J)
\leq T_j(I) \leq 2j$, $T_j(I) = 2j, 1 \leq j \leq m-1$, contradicting the
assumption that $\delta_{m-1}(I) = 1$; hence $T_{m-2}(J) = T_{m-2}(I)$.
We may assume that $T_{m-2}(J) = \deg \mathrm{lcm} (f_2, \cdots, f_{m-1})$.
Then $T_{m-1}((f_2, \cdots, f_m)) = 2 + T_{m-2}(J) > T_{m-1}(J) =
T_{m-1}(I)$, leading to a contradiction.

From the above discussion, and since $T_1 = 2$, clearly there exists $\rho$
such that 
\[
T_l = 
\begin{cases}
2l, & 1 \leq l \leq \rho \\
\min\{\rho + l, n\}& \rho \leq l \leq c
\end{cases}
\]
What we need to show is that $\rho$ is the length of the longest
$R$-regular sequence in $\{f_1, \cdots, f_m\}$. If $f_{j_1}, \cdots,
f_{j_t}$ form a regular sequence, then $T_t = 2t$, so $\rho \geq t$.
Conversely, since $T_\rho = 2\rho$, there exists a regular sequence of
length $\rho$ in $\{f_1, \cdots, f_m\}$. 

Let $y \in V$ be such that $xy \in I$. If $f_1, \cdots, f_s$ are all the
quadratic minimal generators of $(I:x)$, then none of them involves $x$ and
and $y$; therefore, to any regular sequence in $\{f_1, \cdots, f_s\}$, one
can add $xy$, to get a longer regular sequence. The last statement follows
immediately. 
\end{proof}

\begin{lemma}
\label{thm:halfcregseq}
With notation as above, $\rho(I) \geq \frac{c}2$.
\end{lemma}

\begin{proof}
Since $\rho(I) \geq 1$, this holds when $c=1$. By induction on $c$, we may
assume that for all square-free monomial ideals $J$ with $\height J < c$,
$\rho(J) > \frac{\height J}2$. Take a minimal generator $xy$ of $I$. Let $J
= (I \cap \Bbbk[V \setminus \{x,y\}])R$. Since $xy$ is a non-zerodivisor on
$R/J$, $\rho(J) = \rho(I) - 1$, and, further, since, $(J,xy) \subseteq I$,
$\height J < \height (J,xy) \leq \height I$ and  Since $(I,x,y) = (J,x,y)$,
$\height J \geq c-2$. By induction, $\rho(J) \geq \frac{c-2}2$, and,
therefore, $\rho(I) \geq  \frac{c}2$.
\end{proof}

We now prove that Conjecture~\eqref{conj:TB} holds for quadratic monomial
ideals.
{\def\thethm{\ref{thm:quadTB}}
\begin{thm}
\thmquadTB
\end{thm}
\addtocounter{thm}{-1}}

\begin{proof}
We proceed by induction on $c$. If $c=2$, the Taylor bound holds for
$I$~\cite[Corollary~4.3]{HeSrMonBnds04}, so let $c \geq 3$. As discussed in
the previous section, we take $I$ to be the edge ideal of a graph $G$ and
assume that Hypothesis~\ref{standing-hyp} holds. 

For all $x \in V$, notice that $e(R/(I,x))$ is the number of unmixed primes
$\mathfrak p$ of $R/I$ containing $x$.  Since each such prime has height $c$,
in the sum $\sum_{x\in V}{e(R/(I,x))}$, it is counted $c$ times. Therefore 
\[
e(R/I) = \frac1c\sum_{x \in V}{e(R/(I,x)).} 
\] Now suppose $T_c = n$. As noted earlier, $(I,x)$ is essentially of
height $\leq c-1$. Therefore, by induction and by
Proposition~\ref{thm:MlTlcolons}\eqref{thm:TlDecr},
\[
e(R/I) \leq \frac{n}c \frac{T_1T_2 \cdots T_{c-1}}{(c-1)!}
 = \frac{T_1T_2 \cdots T_c}{c!}.
\]
Therefore we may further assume that $T_c = c+ \rho(I) < n$. 

We now reduce to the case that $\rho(I) < c$. If $\rho(I)= c$ then, without loss
of generality, take $f_1, \cdots, f_c$ to be a regular sequence. Write $J =
(f_1, \cdots, f_c)$. Since $J \subseteq I$ and $\height J = c = \height I$,
we see that $e(R/I) \leq e(R/J) = 2^c$. From
Lemma~\ref{thm:first-two-then-one}, $T_l = 2l$ for all $1 \leq l \leq c$.
Hence 
\[
e(R/I) \leq \frac{T_1T_2 \cdots T_c}{c!}.
\]
In particular $G$ is not a collection of $c$ isolated edges, which would
have given $\rho(I) = c$ and $|V| = 2c$. We pick $x \in V$ such that
$\deg_G x \geq 2$. Then $(I:x)$ is essentially of height $\leq c-2$.
Moreover $\rho((I:x)) < \rho(I)$, by Lemma~\ref{thm:first-two-then-one}. We
noted earlier that $(I,x)$ is essentially of height $\leq c-1$. Let $\rho'
\defeq \rho((I,x))$.  Hence, by induction on $c$ and
by Hypothesis~\ref{standing-hyp}, we have
\[
e(R/(I,x)) \leq 
\frac{2 \cdot 4 \cdots 2\rho' \cdot(2\rho'+1)\cdots(c+\rho'-1)}{(c-1)!}
= 2^{\rho'} \mu(\rho',c-1),
\]
which gives, after successive application of~\eqref{eqn:ratios-of-mu-r-c},
(which is permitted since $\rho(I) < c$), $e(R/(I,x)) \leq 2^{\rho(I)}
\mu(\rho(I),c-1)$. Since $\deg_G x \geq 2$ and $\rho((I:x)) \leq
\rho(I)-1$, we can conclude, by a similar argument, that $e(R/(I,x)) \leq
2^{\rho(I)-1} \mu(\rho(I)-1,c-2)$. (Notice that since $\rho(I) - 1 \leq
c-2$, we can apply~\eqref{eqn:ratios-of-mu-r-c}.)

We must show that
\[
e(R/I) \leq 
\frac{2 \cdot 4 \cdots 2\rho(I) \cdot(2\rho(I)+1)\cdots(c+\rho(I))}{c!}
= 2^{\rho(I)} \mu(\rho(I),c).
\]
Since $e(R/I) = e(R/(I,x)) + e(R/(I:x))$, it suffices to show that
\[
2^{\rho(I)} \mu(\rho(I), c-1) + 2^{\rho(I)-1} \mu(\rho(I)-1,c-2) <
2^{\rho(I)}\mu(\rho(I),c).
\]
Set $\rho =
\rho(I)$, $\gamma = c$, $\gamma_1 = 2$. Since $\frac{c}2 \leq \rho(I) < c$,
and $c \geq 3$, we see that $2 \leq \rho < \gamma \leq \rho\gamma_1$ and
$\rho-1 \leq \gamma - \gamma_1$. Applying
Lemma~\ref{thm:mu-r-c-c1-numer-crit} now finishes the proof. 
\end{proof}

\section{Reduction to the Cohen-Macaulay Case of
Theorem~\protect{\ref{thm:quadHHSbipart}}}
\label{sec:rednToCM}

Let $I$ be the edge ideal of a bipartite graph $G$ on $V = V_1 \bigsqcup
V_2$. In order to prove Theorem~\protect{\ref{thm:quadHHSbipart}} for $I$, 
we will first reduce to the case of $G$ having perfect matching, and, by
associating a certain directed graph to $G$, show that only Cohen-Macaulay
bipartite graphs matter. The next section is devoted to proving the theorem
for Cohen-Macaulay bipartite graphs.

A \emph{matching} in $G$ is a maximal (under inclusion) set
$\mathrm m$ of edges such that for all $x \in V$, at most one edge in
$\mathrm m$ is incident on $x$. Edges in a matching form a regular
sequence on $R$. We say that $G$ has \emph{perfect matching}, or, is
\emph{perfectly matched}, if there is a matching $\mathrm m$ such that for
all $x \in V$, there is exactly one edge in $\mathrm m$ is incident on $x$.
K\"onig's theorem~\cite[Section~6.4]{VillmonAlg01} states that the maximum
size of any matching equals the minimum size of any vertex cover.  In the
language of algebra, we can restate it as that the maximum length of a
regular sequence in the set of monomial minimal generators of the edge
ideal equals the height of the ideal. 

\begin{lemma}
\label{thm:GPerfMatchIffCardsOfPartsIsHeight}
With notation as above, $G$ is perfectly matched if and only if  $|V_1| =
|V_2| = \height I = c$.
\end{lemma}

\begin{proof}
If $G$ is perfectly matched, then, first, $ |V_1| = |V_2|$. Secondly, the
matching gives a regular sequence of length $|V_1|$ in $I$, so $\height I
\geq |V_1|$. Since $V_1$ is a vertex cover for $G$, $\height I \leq |V_1|$.
Hence $|V_1| = |V_2| = \height I = c$. Conversely, assume that $|V_1| =
|V_2| = \height I = c$. Since $\height I = c$, $G$ has a minimal vertex
cover of $c$ vertices, and, by K\"onig's
theorem~\cite[Section~6.4]{VillmonAlg01}, a matching of $c$ edges. This
gives a bijection between $V_1$ and $V_2$, so $G$ is perfectly matched.
\end{proof}

\begin{propn}
\label{thm:st-hyp-same-vset-sizes}
Let $I$ be the edge ideal of a bipartite graph $G$ on $V = V_1 \bigsqcup
V_2$. Then Hypothesis~\ref{standing-hyp} holds for $I$ if and only if
$G$ is perfectly matched.
\end{propn}

\begin{proof}
If $G$ is perfectly matched, then let $\mathfrak p \defeq (x : x \in V_1)$
and $\mathfrak q \defeq (x : x \in V_2)$. 
By Lemma~\ref{thm:GPerfMatchIffCardsOfPartsIsHeight}, $\height \mathfrak p
= \mathfrak q = c$. For all $x \in V_1$, $(I,x) \subseteq \mathfrak q$ and
$(I:x) \subseteq \mathfrak q$; the case of $x \in V_2$ is similar. Hence we
see that Hypothesis~\ref{standing-hyp} holds for $I$.

Conversely, assume that $G$ is not perfectly matched. Since $V_1$ and $V_2$
are minimal vertex covers for $G$, we see that $|V_1| \geq c$ and that
$|V_2| \geq c$. In light of
Lemma~\ref{thm:GPerfMatchIffCardsOfPartsIsHeight}, we may assume, without
loss of generality, that $|V_1| > c$. In the paragraph preceding 
Lemma~\ref{thm:GPerfMatchIffCardsOfPartsIsHeight} we noted that there is a
matching with $c$ edges. Let $\{x_1, \cdots, x_c\} \subseteq V_1, \{y_1,
\cdots, y_c\} \subseteq V_2$ be such that $x_1y_1, \cdots, x_cy_c$ is a
matching of $G$. Pick $x \in V_1 \setminus \{x_1, \cdots, x_c\}$. Then $
x_1y_1, \cdots, x_cy_c, x$ is a regular sequence in $(I, x)$, giving
$\height (I, x) > c$. Hence Hypothesis~\ref{standing-hyp} does not hold.
\end{proof}

\begin{remark}
\label{rmk:HtIxIscIsEnough}
The proof above shows that, if $I$ is the edge ideal of a bipartite graph
such that $\height (I, x) = c$ for all $x \in V$, then, $\height (I:x) =
c$, for all $x \in V$. This is not true for arbitrary square-free monomial
ideals.
\end{remark}

\begin{discussion}
\label{disc:defn-of-digraph} For the rest of this section, we restrict our
attention to bipartite graphs $G$ with perfect matching. Let $V_1 = \{x_1,
\cdots, x_c\}$ and $V_2 = \{y_1, \cdots, y_c\}$. We abbreviate $\{1,
\cdots, c\}$ as $[c]$. After relabelling the vertices, we will assume that
$x_iy_i$ is an edge for all $i \in [c]$. We associate $G$ with a directed
graph $\mathfrak d_G$ on $[c]$  defined as follows: for $i, j \in [c]$,
$ij$ is an edge of $\mathfrak d_G$ if and only if $x_iy_j$ is an edge of
$G$.  (Here, by $ij$, we mean the the directed edge from $i$ to $j$.) We
will write $j \succ i$ if there is a directed path from $i$ to $j$ in
$\mathfrak d$. By $j \succcurlyeq i$ we mean that $j \succ i$ or $j = i$.
Let $\mathfrak d$ be any directed graph on $[c]$, and denote the underlying
undirected graph of $\mathfrak d$ by $|\mathfrak d|$.  A vertex $i$ of
$\mathfrak d$ is called a \emph{source} (respectively, \emph{sink}) vertex
if it has no edge directed towards (respectively, away from) it. We say
that a set $A \subseteq [c]$ is an \emph{antichain} if for all $i, j \in
A$, there is no directed path from $i$ to $j$ in $\mathfrak d$, and, by
$\mathcal A_\mathfrak d$, denote the set of antichains in $\mathfrak d$. We
consider $\emptyset$ as an antichain. A \emph{coclique} of $|\mathfrak d|$
is a set $A \subseteq [c]$ such that for all $i \neq j \in A$, $i$ and $j$
are not neighbours in $|\mathfrak d|$. Antichains in $\mathfrak d$ are
cocliques in $|\mathfrak d|$, but the converse is not, in general, true. We
say that $\mathfrak d$ is \emph{acyclic} if there are no directed cycles,
and \emph{transitively closed} if, for all $i, j, k \in [c]$, whenever $ij$
and $jk$ are (directed) edges in $\mathfrak d$, $ik$ is an edge. Observe
that $\mathfrak d$ is a poset under the order $\succ$ if (and only if) it
is acyclic and transitively closed. In this case, for all $A \subseteq
[c]$, $A$ is an antichain in $\mathfrak d$ if and only if $A$ is a coclique
in $|\mathfrak d|$. Let $\kappa(G)$ denote the largest size of any coclique
in $|\mathfrak d_G|$. 
\end{discussion}

Before we proceed, we need the characterization of Cohen-Macaulay bipartite
graphs, due to Herzog-Hibi.
\begin{thm}{\protect{\cite[Lemma~3.3 and Theorem~3.4]{HeHiCMbip05}}}
\label{thm:CMbipHeHi}
Let $G$ be a bipartite graph on $V_1 \bigsqcup V_2$, with edge ideal $I$.
Then $G$ is Cohen-Macaulay if and only if $|V_1| = |V_2| = c = \height I$
and we can write $V_1 = \{x_1, \cdots, x_c\}$ and $V_2 = \{y_1, \cdots,
y_c\}$ such that
\begin{enumerate}
\item \label{thm:CMbipHeHixiyi}
For all $1 \leq i \leq n$, $x_iy_i$ is an edge of $G$.
\item \label{thm:CMbipHeHixiyj}
For all $1 \leq i,j \leq n$, if $x_iy_j$ is an edge of $G$, then $j
\geq i$.
\item \label{thm:CMbipHeHixiyk}
For all $1 \leq i,j,k \leq n$, if $x_iy_j$ and $x_jy_k$ are edges of
$G$, then $x_iy_k$ is an edge of $G$.
\end{enumerate}
\end{thm}

\begin{remark}
\label{rmk:GCMIffDAG}
When we say that $G$ is a Cohen-Macaulay bipartite graph on the vertex set
$\{x_1, \cdots, x_c\} \bigsqcup \{y_1, \cdots, y_c\}$, we will assume that
the variables have already been relabelled so that the conditions of
Theorem~\ref{thm:CMbipHeHi} hold. It is clear that $G$ is Cohen-Macaulay if
and only if $\mathfrak d_G$ is a poset. 
\end{remark}

\begin{lemma}
\label{thm:OnlyXorOnlyYonChains}
Let $G$, $I$, and $\mathfrak d_G$ be as in
Discussion~\ref{disc:defn-of-digraph}. Let $j \succcurlyeq i$. Then for all
$\mathfrak p \in \unm R/I$, if $y_i \in \mathfrak p$, then $y_j \in
\mathfrak p$.
\end{lemma}

\begin{proof}
Applying induction on the length of a directed path from $i$ to $j$, we may
assume, without loss of generality, that $ij$ is a directed edge of
$\mathfrak d_G$. Let $\mathfrak p \in \unm R/I$ and $k \in [c]$. Since
$x_ky_k \in I$, $x_k \in \mathfrak p$ or $y_k \in \mathfrak p$. Since
$\height \mathfrak p = c$, in fact, $x_k \in \mathfrak p$ if and only if
$y_k \not \in \mathfrak p$. Now since $y_i \in \mathfrak p$, $x_i \not \in
\mathfrak p$, so $(I:x_i) \subseteq \mathfrak p$. Note that since $x_iy_j$
is an edge of $G$, $y_j \in (I:x_i)$.
\end{proof}

\begin{discussion}[Collapsing directed graphs]
\label{disc:collapsing-digraphs}
Suppose that $\mathfrak d_G$ has a directed cycle, \textit{i.e.}, a sequence of
directed edges $i_0i_1, i_1i_2, \cdots, i_{p-1}i_p, i_pi_0$, for some $p
\geq 1$. Set $B = \{i_1, \cdots, i_p\}$.  We
\emph{collapse} $\mathfrak d_G$ to obtain a new directed graph
$\tilde{\mathfrak d}$ on the vertex set $[c] \setminus B$ as follows: if
$ij$ an edge of $\mathfrak d_G$ for some $i,j \not \in B$, then $ij$ is an
edge of $\tilde{\mathfrak d}$. For all $1 \leq s \leq p$, and for all $j
\not \in B$, if $i_sj$ (respectively, $ji_s$) is an edge of $\mathfrak
d_G$, then set $i_0j$ (respectively, $ji_0$) to be an edge of
$\tilde{\mathfrak d}$. Any cycle of $\tilde{\mathfrak d}$ comes from a
cycle of $\mathfrak d$; hence the total number of cycles decreases. Let
$\tilde G$ be the bipartite graph associated to $\tilde{\mathfrak d}$. If
$A \subseteq [c] \setminus B$ is a coclique in $|\tilde{\mathfrak d}|$,
then it is clearly a coclique in $|\mathfrak d_G|$; hence $\kappa(\tilde G)
\leq \kappa(G)$. Write $\tilde I \subseteq R$ for the edge ideal of $\tilde
G$. Let $\mathfrak p \in \unm R/I$.  It is an immediate corollary to
Lemma~\ref{thm:OnlyXorOnlyYonChains} that $x_{i_0} \in \mathfrak p$ if and
only if $(x_{i_0}, x_{i_1}, \cdots, x_{i_p}) \subseteq \mathfrak p$, which
holds if and only if $y_{i_s} \not \in \mathfrak p$ for all $0 \leq s \leq
p$. Let $\tilde{\mathfrak p} \defeq (\mathfrak p \cap \Bbbk[\{x_i, y_i : i
\not \in B\}])R$. We first claim that $\tilde I \subseteq \tilde{\mathfrak
p}$. To prove this, we only need to consider the new edges introduced in
$\tilde G$, which are of the form $x_{i_0}y_j$ or $x_jy_{i_0}$ for some $j
\in [c] \setminus B$. For the edge $x_{i_0}y_j$, if $x_{i_0} \not \in
\tilde{\mathfrak p}$, then $x_{i_0} \not \in \mathfrak p$; by the above
observation, we see that $x_{i_s} \not \in \mathfrak p$ for all $0 \leq s
\leq p$. Hence $y_j \in \mathfrak p$ giving $y_j \in \tilde{\mathfrak p}$.
The case of $x_jy_{i_0}$ is similar. Now since $\height \tilde{\mathfrak p}
= \height \tilde I = c - |B|$, we conclude that $\tilde{\mathfrak p} \in
\unm R/{\tilde I}$. The map $\unm R/I \rightarrow \unm R/\tilde I$ sending
$\mathfrak p \mapsto \tilde{\mathfrak p}$ is injective. Conversely, let 
$\tilde{\mathfrak q} \in \unm R/\tilde I$. Set 
\[
\mathfrak q  \defeq 
\begin{cases}
\tilde{\mathfrak q} + (x_{i_1}, \cdots, x_{i_p}) & \text{if}\; x_{i_0} \in
\tilde{\mathfrak q} \\
\tilde{\mathfrak q} + (y_{i_1}, \cdots, y_{i_p}) & \text{if}\; y_{i_0} \in
\tilde{\mathfrak q}.
\end{cases}
\]
Then we get an injective map $\unm R/\tilde I \rightarrow \unm R/I$.
Therefore, $e(R/I) = e(R/\tilde I)$. 
\end{discussion}

\begin{discussion}[Closing directed graphs under transitivity]
\label{disc:transitively-closing}
Suppose that $ij$ and $jk$ are edges of $\mathfrak d_G$; then we add an
edge $ik$. Call the new graph $\widehat{\mathfrak d}$ and let $\widehat G$
be the bipartite graph associated to $\widehat{\mathfrak d}$. Let $\widehat
I$ be the edge ideal of $\widehat G$. Since $I \subseteq \widehat I$ and
$\height I = \height \widehat I$, we have that $e(R/I) \geq e(R/\widehat
I)$. In order to show that $e(R/I) = e(R/\widehat I)$, it suffices to show
that $x_iy_k \in \mathfrak p$, for all $\mathfrak p \in \unm R/I$. Let
$\mathfrak p \in \unm R/I$ be such that $x_i \not \in \mathfrak p$. Then,
since $k \succ i$, by Lemma~\ref{thm:OnlyXorOnlyYonChains}, $y_k \in
\mathfrak p$, and therefore, $x_iy_k \in \mathfrak p$. Moreover, any
coclique in $|\widehat{\mathfrak d}|$ is a coclique in $|\mathfrak d_G|$,
so $\kappa(\widehat G) \leq \kappa(G)$.
\end{discussion}

The significance of $\kappa(G)$ is that it gives a lower bound on the
(Castelnuovo-Mumford) regularity, $\reg R/I$. Following
Zheng~\cite{Zhereslnfacets04}, we say that two edges $vw$ and $v'w'$ of a
graph $G$ are \emph{disconnected} if they are no more edges between the
four vertices $v,v',w,w'$. The edges in any pairwise disconnected set form
a regular sequence in $R$; in fact, a set $\mathbf a$ of edges is pairwise
disconnected if and only if $(I \cap \Bbbk[V_\mathbf a])R$ is generated by
the regular sequence of edges in $\mathbf a$, where by $V_\mathbf a$, we
mean the set of vertices on which the edges in $\mathbf a$ are incident.
The latter condition holds if and only if the subgraph of $G$ induced on
$V_\mathbf a$, denoted as $G|_{V_\mathbf a}$, is a collection of $|\mathbf
a|$ isolated edges. Set $r(I) \defeq \max \{ |\mathbf a| : \mathbf a \;
\text{is a set of pairwise disconnected edges in} \; G\}$. If $G$ is a
forest with edge ideal $I$, then $\reg R/I =
r(I)$~\cite[Theorem~2.18]{Zhereslnfacets04}. 

\begin{lemma}
\label{thm:rIIsAtLeastKappaG}
With notation as above, $r(I) \geq \kappa(G) \geq \max \{|A| : A \in
\mathcal A_{\mathfrak{d}_G}\}$. 
\end{lemma}

\begin{proof}
If $A \subseteq [c]$ is a coclique of $|\mathfrak d_G|$, we easily see that
the edges $\{x_iy_i : i \in A\}$ are pairwise disconnected in $G$. The
assertion now follows from the observation, that we made in
Discussion~\ref{disc:defn-of-digraph}, that any antichain in $\mathfrak
d_G$ is a coclique of $|\mathfrak d_G|$.
\end{proof}

\begin{lemma}
\label{thm:Ml-geq-r-plus-l}
With notation as above, for $1 \leq l \leq r(I)$, $M_l(I) = 2l$ and for $r(I)
\leq l \leq c$, $M_l(I) \geq l + r(I)$. Hence, for all $1 \leq l \leq c$,
$M_l(I) > l$.
\end{lemma}

\begin{proof}
Let $\mathbf a$ be a set of pairwise disconnected edges with $|\mathbf a| =
r(I)$. Then, with the notation as above, $(I\cap\Bbbk[V_{\mathbf a}])R$ is
generated by a regular sequence of length $r(I)$. From
Lemma~\ref{thm:varSetRestrictColon}\eqref{thm:varSetRestrict}, we have
that, for $1 \leq l \leq r(I)$, $M_l(I) \geq 2l$. From the the Taylor
resolution of
$R/I$, it follows that $M_l \leq 2l$. Hence $M_l = 2l$ for all $1 \leq l
\leq r(I)$.  For $l > r(I)$, we see from Lemma~\ref{thm:MlincrHeZh} that
$M_l \geq l+r(I)$.
\end{proof}

\begin{discussion}[Reduction to the Cohen-Macaulay case]
\label{disc:redn-to-cm-bip}
Now let $G$ be any perfectly matched bipartite graph. We first collapse
$\mathfrak d_G$, repeatedly if necessary, to get a directed acyclic graph,
which we denote $\tilde{\mathfrak d}$. We now close $\tilde{\mathfrak d}$
under transitivity, and call it $\widehat{\mathfrak d}$. Denote the
corresponding bipartite graph by $\widehat G$, and its edge ideal by
$\widehat I$. Notice that $\widehat G$ is Cohen-Macaulay, from
Remark~\ref{rmk:GCMIffDAG}. From the discussion, we see that $\height
\widehat I \leq c$, $e(R/\widehat I) = e(R/I)$ and that $\kappa(\widehat G)
\leq \kappa(G)$.
\end{discussion}

Here is the outline of the rest of the proof: since
$\widehat G$ is Cohen-Macaulay, $r(\widehat I) = \kappa(\widehat G)$ and
that equality must also hold for $\widehat I$ in
Lemma~\ref{thm:Ml-geq-r-plus-l}.  Hence $M_l(\widehat I) \leq M_l(I)$ for
$1 \leq l \leq \height \widehat I$.  Now, since $\height \widehat I < l
\leq c$, $M_l(I) > l$, the conjectured bound for $I$ would be established,
if it can be established for $\widehat I$.

\section{Cohen-Macaulay Bipartite Graphs}
\label{sec:cm-bipartite-resln}

Lemma~\ref{thm:Ml-geq-r-plus-l} gives that $\reg R/I \geq r(I)$.
For arbitrary bipartite graphs, this might be a strict inequality,
(consider, \textit{e.g.}, the edge ideal of the cycle on $8$ vertices), but
we have:
\begin{propn}
\label{thm:regIsrIforCM}
Let $I$ be the edge ideal of be a Cohen-Macaulay bipartite graph $G$ on the
vertex set $\{x_1, \cdots, x_c\} \bigsqcup \{y_1, \cdots, y_c\}$. Then
$\reg R/I = r(I)$.
\end{propn}

\begin{proof}
It suffices to show that $\reg R/I \leq r(I)$, by induction on the number
of vertices. Since the claim is true for a Cohen-Macaulay bipartite graph
on $2$ vertices, we assume inductively that for all Cohen-Macaulay
bipartite graphs on fewer than $2c$ vertices, the claim holds.  

Note that $y_1$ is a leaf vertex of $G$. Since $\depth
\frac{R}{(I:x_1)} \geq \depth R/I = \dim R/I = \dim \frac{R}{(I:x_1)}$
(Remark~\ref{rmk:ColonHasSameDepth} and Hypothesis~\ref{standing-hyp}),
$\frac{R}{(I:x_1)}$ is Cohen-Macaulay. Further, $(I,x_1)$
is the edge ideal of the deletion of the vertices $x_1$ and $y_1$ in $G$;
this graph satisfies the conditions in
Theorem~\ref{thm:CMbipHeHi}, so $\frac{R}{(I,x_1)}$ is Cohen-Macaulay.
Moreover $\height (I:x_1) = c = \height (I,x_1)$. Let $\mathbf a$ be a set
of pairwise disconnected edges in the graph of $(I:x_1)$. Then $\mathbf a$
is pairwise disconnected in $G$, because the graph $(I:x_1)$ is obtained by
deleting the neighbours of $x_1$ from $G$.  Hence we may join $x_1y_1$ to
get a set of pairwise disconnected edges in $G$, so $r((I:x_1)) \leq r(I) -
1$.  Let $J = (I \cap \Bbbk[x_2, \cdots, x_c, y_2, \cdots, x_c])R$. Then
$J$ is the edge ideal of the deletion $G \setminus x_1$ and $(I,x_1) =
(J,x_1)$.  It is evident that $r((I,x_1)) = r(J) \leq r(I)$. 

We need to show that $M_l(I) \leq l + r(I)$ for all $1 \leq l \leq c =
\height I = \projdim R/I$. From the exact sequence
\[
\xymatrix{%
\ar[r] & \tor_l\left(\Bbbk,\frac{R}{(I:x_1)}(-1)\right) \ar[r] & 
\tor_l(\Bbbk,R/I) \ar[r] & 
\tor_l\left(\Bbbk,\frac{R}{(I,x_1)}\right) \ar[r] &,}
\]
we can see that showing
\[
\tor_l\left(\Bbbk,\frac{R}{(I:x_1)}(-1)\right)_j =0=
\tor_l\left(\Bbbk,\frac{R}{(I,x_1)}\right)_j 
\; \text{for all}\; j > l + r(I),
\]
will suffice. This is equivalent, by the induction hypothesis, to showing
that $r((I:x_1)) \leq r(I) - 1$ and that $r((I,x_1) \leq r(I)$, which we
have done.
\end{proof}

\begin{cor}
\label{thm:Ml-is-r-plus-l-CM}
With notation as above, for $1 \leq l \leq r(I)$, $M_l(I) = 2l$ and for $r(I)
\leq l \leq c$, $M_l(I) = l + r(I)$.
\end{cor}

\begin{proof}
Follows from Lemma~\ref{thm:Ml-geq-r-plus-l} and the definition of
regularity.
\end{proof}

For the rest of this section, we will take $I$ to be the edge ideal of 
an arbitrary Cohen-Macaulay bipartite graph $G$ on the vertex set $\{x_1,
\cdots, x_c\} \bigsqcup \{y_1, \cdots, y_c\}$. Recall that the labelling of
vertices was chosen so that the conditions of Theorem~\ref{thm:CMbipHeHi}
hold; see Remark~\ref{rmk:GCMIffDAG}. Hence for all $i, j \in [c]$, if $j
\succ i$ then $j > i$.

\begin{propn}
\label{thm:rI-is-max-antich-size}
With notation as above, $r(I) = \max\{|A| : A \in \mathcal
A_{\mathfrak{d}_G}\}$. Consequently, $r(I) = \kappa(G)$.
\end{propn}

\begin{proof}
Let $A \in \mathcal A_{\mathfrak{d}_G}$. Then $\{x_iy_i: i \in A \}$ is a
set of pairwise disconnected edges in $G$. Conversely, let $\mathbf b$ be a
set of pairwise disconnected edges such that there exists $j \neq i$ such
that $x_iy_j \in \mathbf b$. Let $\mathbf a \defeq \left(\mathbf b
\setminus \{x_iy_j\}\right) \cup \{x_iy_i\}$. We claim that the edges in
$\mathbf a$ are pairwise disconnected; for, if not, then some edge in 
$\mathbf b \setminus \{x_iy_j\}$ is incident on a neighbour of $x_i$ or
$y_i$. We claim that this must be on a neighbour of $y_i$, for, if it were
on a neighbour of $x_i$, then the set $\mathbf b$ would not have been
pairwise disconnected. Therefore $k < i$ such that $x_ky_i$ is an edge of
$G$, and some edge incident on $x_k$ belongs to $\mathbf a$.  However,
since $G$ is Cohen-Macaulay, $x_ky_j$ is an edge of $G$ too, contradicting
the hypothesis that the edges of $\mathbf b$ are pairwise disconnected.
Repeating this if necessary, we can construct a set $\mathbf a$ of pairwise
disconnected edges in $G$ such that $\mathbf a = \{x_iy_i : i \in A\}$ for
some $A \subseteq [c]$ and $|\mathbf a| = |\mathbf b|$. Such a set $A$ is
an antichain in $\mathfrak d_G$.
\end{proof}

\begin{propn}
\label{thm:multIsAsManyAntiCh}
With notation as above, $e(R/I) = |\mathcal A_{\mathfrak{d}_G}|$.
\end{propn}

\begin{proof}
Let $\mathfrak p \in \unm R/I$. Let $A \defeq \{i \in [c] : y_i \in
\mathfrak p \; \text{and for all} \; j \in [c] \; \text{with} \; i \succ j
, y_{j} \not \in \mathfrak p\}$. Note that $A$ is an antichain. This gives
a map from $\unm R/I$ to $\mathcal A_{\mathfrak{d}_G}$, which is injective
by Lemma~\ref{thm:OnlyXorOnlyYonChains}. Conversely, for any antichain $A$
of $\mathfrak d_G$, the prime ideal $(x_j : j \not \succcurlyeq i\;
\text{for any} \; i \in A)  + (y_j : j \succcurlyeq i\; \text{for some} \;
i \in A)$ belongs to $\unm R/I$. This gives a bijection $\mathcal
A_{\mathfrak{d}_G}$ and $\unm R/I$, with the empty set corresponding to
$(x_1, \cdots, x_c)$.
\end{proof}

\begin{lemma}
\label{thm:ManyAntiChBd}
Let $\mathfrak d$ be any poset on $c$ vertices, with order $\succ$,
$\mathcal A$ the set of antichains in $\mathfrak d$ and $r = \max\{|A| : A
\in \mathcal A\}$. Then $|\mathcal A| \leq 2^r \mu(r,c)$. Equality holds
above, if and only if $r=1$ or $r=c$.
\end{lemma}

\begin{proof}
We prove this by induction on $c$. If $r=1$, (in particular, if $c=1$),
$\mathfrak d$ is a chain, \textit{i.e.}, for all $i \neq j \in [c]$, $i
\succ j$ or $j \succ i$. In this case, $|\mathcal A| = c+1 = 2\mu(1,c)$.
If $c=r\geq2$, then $\mathfrak d$ is a collection of $c$ isolated vertices,
in which every subset of $[c]$ is an antichain, \textit{i.e.}, $|\mathcal
A| = 2^c = 2^c\mu(c,c)$. Note that equality holds in both the cases above.

We now have $c > r \geq 2$. Pick a vertex $i$ such that there is an
antichain $A$ with $i \in A$ and $|A| = r$. Set $\tilde{\mathfrak d} \defeq
\{j \in \mathfrak d: j \not \succcurlyeq i \, \text{or} \, i \not
\succcurlyeq j\}$. Let $\mathfrak d'$ be the poset obtained by deleting $i$
from $\mathfrak d$, keeping all the other elements and relations among
them. Denote the respective sets of antichains by $\tilde{\mathcal A}$ and
$\mathcal A'$. Now for any $A \subseteq [c]$, $A \in \mathcal A \setminus
\mathcal A'$ if and only if $i \in A$ and $A \setminus \{i\} \in
\tilde{\mathcal A}$. Therefore $\mathcal A = \mathcal A' \bigsqcup \{A \cup
\{i\} : A \in \tilde{\mathcal A}\}$ and $|\mathcal A| = |\mathcal A'| +
|\tilde{\mathcal A}|$. 

Observe that $\max\{|A| : A \in \tilde{\mathcal A}\} = r-1$.
Let $r' \defeq \max\{|A| : A \in \mathcal A'\}$. Then $r'
\leq r$. Let $c_1 \defeq | \{ j \in \mathfrak d : j \succcurlyeq i \,
\text{or} \, i \succcurlyeq j\}|$. Then $\tilde{\mathfrak d}$ has $c-c_1$
vertices. We note that $r-1 \leq c-c_1$. We assume, by induction on the
number of vertices, that the lemma holds for $\tilde{\mathfrak d}$ and
$\mathfrak d'$, yielding
\[
|\mathcal A| \leq 2^{r'}\mu(r',c-1)+2^{r-1}\mu(r-1,c-c_1),
\]
and, by repeated application of~\eqref{eqn:ratios-of-mu-r-c} from
Discussion~\ref{disc:defn-of-mu-r-c}, (which is permitted since $r' \leq r
\leq c-1$)
\begin{equation}
\label{eq:no-of-antich}
|\mathcal A|  \leq 2^r\mu(r,c-1)+2^{r-1}\mu(r-1,c-c_1).
\end{equation}
Since $c > r \geq 2$, we must show that $|\mathcal A_{\mathfrak d}| <
2^r\mu(r,c)$; to this end, it suffices to show that
\[
2^r\mu(r,c-1)+2^{r-1}\mu(r-1,c-c_1) < 2^r\mu(r,c),
\]
which follows from Lemma~\ref{thm:mu-r-c-c1-numer-crit} with $\rho = r,
\gamma = c, \gamma_1 = c_1$. Note that by the choice of $i$, $c \leq rc_1$.
\end{proof}

{\def\thethm{\ref{thm:quadHHSbipart}}
\begin{thm}
\thmquadHHSbipart
\end{thm}
\addtocounter{thm}{-1}}

\begin{proof}
Denote the bipartite graph by $G$ and its vertex set by $V = V_1 \bigsqcup
V_2$. We may assume that Hypothesis~\ref{standing-hyp} holds;
hence $G$ is perfectly matched. We reduce the proof to the Cohen-Macaulay
case to obtain $\widehat G$ and $\widehat I$ as in
Discussion~\ref{disc:redn-to-cm-bip}. Let $\widehat c \defeq \height
\widehat I$. The height of the edge ideal does not increase during
collapsing the directed graph, while it remains unchanged after closing the
directed graph under transitivity, so $\widehat c \leq c$. Again, from
Discussion~\ref{disc:redn-to-cm-bip}, $e(R/I) = e(R/\widehat I)$. We
already observed that $r(I) \geq \kappa(G) \geq \kappa(\widehat G) =
r(\widehat I)$; see Lemma~\ref{thm:rI-is-max-antich-size}. Now from
Lemma~\ref{thm:Ml-geq-r-plus-l} and Corollary~\ref{thm:Ml-is-r-plus-l-CM},
we see that, for $1 \leq l \leq \widehat c$, $M_l(I) \geq M_l(\widehat I)$
and that $M_l(I) > l$ for all $\widehat c < l \leq c$. Hence it suffices to
show that 
\[
e(R/\widehat I ) \leq \frac{M_1(\widehat I) \cdots M_{\widehat
c}(\widehat I)}{\widehat c!}.
\]
From Proposition~\ref{thm:multIsAsManyAntiCh}, $e(R/\widehat I) = |\mathcal
A_{\mathfrak d_{\widehat G}}|$. Corollary~\ref{thm:Ml-is-r-plus-l-CM} gives
\[
\frac{M_1(\widehat I) \cdots M_{\widehat c}(\widehat I)}{\widehat c!} =
2^{r(\widehat I)} \mu(r(\widehat I), \widehat c).
\]
Since, by Proposition~\ref{thm:rI-is-max-antich-size}, $r(\widehat I) = 
\max\{|A| : A \in \mathcal A_{\mathfrak d_{\widehat G}}\}$, we apply
Lemma~\ref{thm:ManyAntiChBd}, with $\mathcal A = \mathcal
A_{\mathfrak d_{\widehat G}}$, to finish the proof.
\end{proof}

\section{Pure and Quasi-pure Resolutions}
\label{sec:quasiPureRes}

When can equality hold for $I$ in the conjectured bound? The proof
Theorem~\ref{thm:quadHHSbipart} above and Lemma~\ref{thm:ManyAntiChBd},
show that if $G$ is a Cohen-Macaulay bipartite graph with edge ideal $I$,
and equality holds for $I$, then $\reg R/I = c$ or $\reg R/I = 1$.
We are now ready to prove Theorem~\ref{thm:quadHHSe}.

{\def\thethm{\ref{thm:quadHHSe}}
\begin{thm}
\thmquadHHSebipart
\end{thm}
\addtocounter{thm}{-1}}

\begin{proof}
Denote the bipartite graph by $G$.  We first reduce to the case that
Hypothesis~\ref{standing-hyp} holds. We will show that $\height (I,x) = c$
for $x \in V$; this suffices, by Remark~\ref{rmk:HtIxIscIsEnough}. 

Assume, by way of contradiction, that $x \in V$ is such that $\height (I,x)
> c$. Then $\height (I:x) = c$ and $e(R/(I:x)) = e(R/I)$. We may assume
that $x$ is not an isolated vertex of $G$; for otherwise, $x$ would not
have divided any minimal generator of $I$. Hence $x$ has at least one
neighbour, so $(I:x)$ is essentially of height at most $c-1$; see the
paragraph following Proposition~\ref{thm:MlTlcolons}. Let $J \subseteq R$
be the ideal generated by the quadratic minimal generators of $(I:x)$.
Observe that $(I:x)$ is generated by the neighbours of $x$, modulo $J$.
Hence $e(R/(I:x)) = e(R/J)$. It follows from
Lemma~\ref{thm:varSetRestrictColon}\eqref{thm:varSetRestrict} and
Proposition~\ref{thm:MlTlcolons}\eqref{thm:MlDecr} that $M_l(J) \leq
M_l((I:x)) \leq M_l(I)$ for all $1 \leq l \leq c$. Now, $M_l(J) > l$, for
all $l$.  Therefore equality holds for $J$ in Conjecture~\eqref{conj:HHSu}.
Since $M_l((I:x)) \geq M_l(J)$ and $\height J < c = \height (I:x)$, we see
that equality cannot hold for $(I:x)$, and hence, again by
Proposition~\ref{thm:MlTlcolons}\eqref{thm:MlDecr}, for $I$. Therefore we
may assume that Hypothesis~\ref{standing-hyp} holds.

By Proposition~\ref{thm:st-hyp-same-vset-sizes}, $G$ has perfect matching.
Let $\mathfrak d_G$ be the directed graph associated to $G$, as in
Discussion~\ref{disc:defn-of-digraph}. We can now reduce the problem to the
Cohen-Macaulay case. Let $\widehat G$ and $\widehat I$ be as in the proof
of Theorem~\ref{thm:quadHHSbipart}. Since equality holds for $I$, we see
that $\widehat c = \height I = c$, $M_l(\widehat I) = M_l(I)$ for all $1
\leq l \leq c$. In particular, since $\reg R/I = \reg R/\widehat I$, it
follows from Proposition~\ref{thm:rI-is-max-antich-size} and
Discussion~\ref{disc:redn-to-cm-bip} that $\kappa(\widehat G) = \kappa(G)$.
Moreover, equality must hold for $\widehat I$.  

Since $\widehat G$ is Cohen-Macaulay, $\mathcal A_{\mathfrak d_{\widehat
G}}$ is a poset, and, from Lemma~\ref{thm:ManyAntiChBd}, we see that
$r(\widehat I)=1$ or $r(\widehat I)=\widehat c$. If $r(\widehat I) = c$, then
$\mathcal A_{\mathfrak d_{\widehat G}} = \mathcal A_{\mathfrak d_G}$ is a
collection of $c$ isolated vertices, or, equivalently, $G$ is a collection
of $c$ isolated edges. In this case, $R/I$ is a complete intersection.
Since all the minimal generators of $I$ have the same degree, $R/I$ has a
pure resolution.

If $r(\widehat I) = 1$, then $\mathcal A_{\mathfrak d_{\widehat G}}$ and,
hence, $\mathcal A_{\mathfrak d_G}$ have precisely one source vertex and
one sink vertex. With that, $1 = r(\widehat I) = \kappa(\widehat G) =
\kappa(G)$ if and only if $\mathcal A_{\mathfrak d_G}$ is a chain,
\textit{i.e.}, $\mathcal A_{\mathfrak d_{\widehat G}} = \mathcal
A_{\mathfrak d_G}$. In other words, $R/I$ is Cohen-Macaulay with $\reg R/I
= 1$, which, evidently, has a pure resolution.
\end{proof}

A Cohen-Macaulay bipartite graph $G$ on the vertex set $\{x_1, \cdots,
x_c\} \bigsqcup \{y_1, \cdots, y_c\}$ has $r(I) = 1$, if and only if, in
the context of Theorem~\ref{thm:CMbipHeHi}, $x_iy_j$ is an edge for all $j
\geq i$. Such a graph $G$ is acyclic (\textit{i.e.}, a forest), if and only
if it is a path on three edges.

Suppose that $G$ is a graph on $V$ with edge ideal $I$ with the property
that $r(I) = \reg R/I$. For instance, $G$ is a forest
(\cite[Theorem~2.18]{Zhereslnfacets04}) or a Cohen-Macaulay bipartite
graph (Proposition~\ref{thm:regIsrIforCM}). If $r(I) \leq 2$, then $R/I$
has a quasi-pure resolution.  Now suppose that $r(I) \geq 3$ and that $R/I$
has a quasi-pure resolution.  Since $r(I) \geq 3$, $M_3(I) = 6$, from
Lemma~\ref{thm:Ml-geq-r-plus-l}.  Hence $m_4(I) \geq 6$. We claim that
every vertex has at most $3$ neighbours. More generally,
\begin{propn}
\label{thm:ml-is-deg}
Let $I$ be the edge ideal of a graph $G$. For any multidegree $\sigma$,
$\beta_{|\sigma|-1, \sigma}(R/I) \neq 0$ if and only if there exists a
partition $\sigma = \sigma_1 \bigsqcup \cdots \bigsqcup \sigma_d$ such that
for all $1 \leq i \leq d$, for all $x \in \sigma_i$ and for all $y \in
\sigma \setminus \sigma_i$, $xy$ is an edge of $G|_\sigma$.
\end{propn}

\begin{proof}
We immediately reduce the problem to the case that $\sigma = V$, noting
that, by Lemma~\ref{thm:varSetRestrictColon}\eqref{thm:varSetRestrict} and
the flatness $R$ over $\Bbbk[\sigma]$,
$\beta_{|\sigma|-1, \sigma}(R/I) = \beta_{|\sigma|-1, \sigma} \left(
\frac{R}{(I \cap \Bbbk[\sigma])R} \right) = \beta_{|\sigma|-1, \sigma}
\left( \frac{\Bbbk[\sigma]}{(I \cap \Bbbk[\sigma])} \right)$ and that $(I
\cap \Bbbk[\sigma])$ is the edge ideal of $G|_\sigma$ in the ring
$\Bbbk[\sigma]$. Hence we need to show that $\beta_{|V|-1, V}(R/I) \neq 0$
if and only if there is a partition $V = V_1 \bigsqcup \cdots \bigsqcup
V_d$ such that for all $1 \leq i \leq d$, for all $x \in V_i$ and for all
$y \in V \setminus V_i$, $xy$ is an edge of $G$. Let $\Delta$ be the
coclique complex of $G$. Hochster's formula gives that $\beta_{|V|-1,
V}(R/I) = \dim_\Bbbk \rhomo_0(\Delta; \Bbbk)$. Hence we must show that
$\Delta$ is disconnected if and only if a partition, such as above,
exists. 

Suppose such a partition exists. Then any coclique of $G$ is contained in
$V_i$ for some $1 \leq i \leq d$; hence $\Delta$ is disconnected.
Conversely, assume that $\Delta$ is disconnected. Denoting the number of
distinct components of $\Delta$ by $d$, we set $V_i, 1 \leq i \leq d$, be
the vertex sets of these components. We see immediately that for $x, y \in
V$, whenever $x$ and $y$ are in different components of $\Delta$, there is
an edge $xy$ in $G$.
\end{proof}

We wish to mention here that this agrees with the result of
Novik-Swartz~\cite[Theorem~1.3]{NoSwCMconn06} that the first skip in the
sequence of $m_l$'s is at $n - q_1 + 1$, where $q_1$ is the Cohen-Macaulay
connectivity of the $1$-dimensional skeleton of the Stanley-Reisner complex
of $I$.  For the edge ideal of a graph $G$, the $1$-dimensional skeleton of
its Stanley-Reisner complex is the complement graph $\bar G$. In passing,
let us note that if $G$ is a forest, then Proposition~\ref{thm:ml-is-deg}
implies that $\max \{l : m_l(I) = l+1\} = \max \{\deg_G x : x \in V\}$.
More generally, if $G$ is a bipartite graph, then $\max \{l : m_l(I) =
l+1\}$ is the largest cardinality of a complete bipartite subgraph of
$G$. 

\begin{discussion}
\label{disc:FirstSkipInml}
We already observed that if $\reg R/I \leq 2$, then $R/I$ has a quasi-pure
resolution. Let $G$ be a connected Cohen-Macaulay bipartite graph such that
$\reg R/I \geq 3$ and $R/I$ has a quasi-pure resolution. It is easy to see
that $G$ is connected if and only if $\mathfrak d_G$ is a connected poset.
Since $\reg R/I = r(I) \geq 3$, $M_3(I) = 6$ by
Corollary~\ref{thm:Ml-is-r-plus-l-CM}, so for $R/I$ to have a quasi-pure
resolution, we must have $m_4(I) \geq 6$. This means, by the observation in
the last paragraph, that for all $i$, there are at most two elements $j$
such that $j \succcurlyeq i$ (or $i \succcurlyeq j$) in $\mathfrak d_G$.
For $i, j \in [c]$, say that $j$ \emph{covers} $i$ if $j \succ i$ and there
does not exist $j'$ such that $j \succneqq j' \succneqq i$. Since
$\mathfrak d_G$ is connected, in every maximal chain, there exists $i, j,
j'$ such that $j$ and $j'$ cover $i$ or $i$ covers $j$ and $j'$. From the
observation above, it follows that, in the first case, $i$ is a source
vertex and that $j$ and $j'$ are sink vertices. Similarly, in the second
case, $i$ is a sink vertex and $j$ and $j'$ are source vertices. Hence
every maximal chain of $\mathfrak d_G$ has length at most one; in fact,
since $\mathfrak d_G$ is connected, every maximal chain has length one.
Therefore every vertex in $\mathfrak d_G$ is a source vertex or a sink
vertex, but not both. Every source (respectively, sink) vertex in
$\mathfrak d_G$ is covered by (respectively, covers) at most two sink
(respectively, source) vertices. For $x_iy_i$ to be a leaf in $G$, it is
necessary and sufficient that $i$ is a source vertex or a sink vertex in
$\mathfrak d_G$. Therefore, in our case, $x_iy_i$ is a leaf for all $i$; in
other words, $G$ is the suspension\footnotemark of its subgraph $G'$
induced on the set of vertices $V' \defeq \{x_i : i\; \text{is a source
vertex of}\; \mathfrak d_G\} \cup \{y_i : i\; \text{is a sink vertex of}\;
\mathfrak d_G\} \subseteq V$. The underlying undirected graph $|\mathfrak
d_G|$ is, first, bipartite, and, secondly, a path on $c$ vertices
(necessarily, if $c$ is odd) or a cycle on $c$ vertices. The subgraph $G'$
of $G$ described above is a path (respectively, a cycle) if and only if
$|\mathfrak d_G|$ is a path (respectively, a cycle).
\end{discussion}

\footnotetext{%
A graph $G$ is said to be the \emph{suspension} of a subgraph $G'$, if $G$
is obtained by attaching exactly one leaf vertex to every vertex of $G'$.
If $G$ is the suspension of a subgraph $G'$, then $G$ is Cohen-Macaulay.
We see this as follows. Let $V' \subseteq V$ be the set of vertices of
$G'$. Denote the edge ideal of $G'$ in $\Bbbk[V']$ by $I'$. Then $I$ is the
polarization of $I' + (x^2 : x \in V')$ in the ring $R = \Bbbk[V]$. Since 
$\Bbbk[V'] /(I' + (x^2 : x \in V'))$ is Artinian, we see that $R/I$ is
Cohen-Macaulay. Villarreal showed that a tree $G$ is Cohen-Macaulay if and
only if $G$ is the suspension of a subgraph $G'$; see,
\textit{e.g.},~\cite[Theorem~6.5.1]{VillmonAlg01}.}

\begin{propn}
\label{thm:qpureCMbip}
Let $I$ be the edge ideal of a connected Cohen-Macaulay bipartite graph
$G$. Then the following are equivalent:
\begin{enumerate}
\item \label{qpureCMbipRegGeqThree} $\reg R/I \geq 3$ and $R/I$ has a
quasi-pure resolution.
\item \label{qpureCMbipAtMostSixVert} $G$ is the suspension of the path on
five or six vertices or of the $6$-cycle.
\end{enumerate}
\end{propn}

\begin{proof}
\eqref{qpureCMbipRegGeqThree} $\implies$ \eqref{qpureCMbipAtMostSixVert}:
First, if $c \geq 7$, then we claim that $R/I$ cannot have a quasi-pure
resolution. Since $\reg R/I \geq 3$, $\mathfrak d_G$ is such that every
vertex is a source vertex or a sink vertex, but not both, and that every
source (respectively, sink) vertex in $\mathfrak d_G$ is covered by
(respectively, covers) at most two sink (respectively, source) vertices. 
If $c>7$, then restrict $\mathfrak d_G$ to one of its connected subgraphs
with seven vertices. This corresponds to restricting $G$ to a
Cohen-Macaulay subtree on $14$ vertices. If we show that the edge ideal of
this subtree does not have a quasi-pure resolution, then, by
Lemma~\ref{thm:varSetRestrictColon}\eqref{thm:varSetRestrict}, we have that
$R/I$ does not have a quasi-pure resolution. Therefore, by replacing $G$ by
this subgraph, we may assume that $G$ is a Cohen-Macaulay tree on $14$
vertices, such that the length of every maximal path in $\mathfrak d_G$ is
one, and prove the $R/I$ does not have a quasi-pure resolution. We may
verify this with a computer algebra system, such as~\cite{M2}, but 
we give a direct proof below.

We will prove this when $\mathfrak d_G$ has four source vertices and three
sink vertices. The other case is of $\mathfrak d_G$ with three source
vertices and four sink vertices; this corresponds to relabelling the
partition of the vertex set. Since $c=7$ is odd, $G$ is the suspension of a
path on $7$ vertices. We label the source vertices $1, \cdots, 4$ and the
sink vertices $5, 6, 7$. Then the edges of $\mathfrak d_G$ are $15, 25, 26,
36, 37$ and $47$. Hence $I = (x_1y_1, \cdots, x_7y_7, x_1y_5, x_2y_5,
x_2y_6, x_3y_6, x_3y_7, x_4y_7)$. We saw that $m_4(I) = 6$. Since the set
of four source vertices in $\mathfrak d_G$ form an antichain,
$\reg R/I = 4$, and hence $M_4(I) = 8$; to prove that $R/I$ does not have a
quasi-pure resolution, we need to show that $m_5(I) \leq 7$. Let $\sigma =
\{x_1, y_5, x_2, y_6, x_3, y_7, x_5\}$, and $J = (I \cap \Bbbk[\sigma])R =
(x_1y_5,x_2y_5,x_5y_5,x_2y_6,x_3y_6,x_3y_7)$. We will show that
$\beta_{5,7}(R/J) \neq 0$, which will suffice, by
Lemma~\ref{thm:varSetRestrictColon}\eqref{thm:varSetRestrict}, to show that 
$m_5(I) \leq 7$. We have a short exact sequence of graded $R$-modules
\[
\xymatrix{%
0 \ar[r] & \frac{R}{(J:y_5)} (-1)) \ar[r] & R/J \ar[r] & R/(J,y_5) \ar[r] &
0.}
\]
Since $R/(J,y_5)$ is Cohen-Macaulay and $\height (J, y_5) = 3$, we see from
the associated long exact sequence of $\tor(\Bbbk,-)$ that
\[
\tor_5(\Bbbk,\frac{R}{(J:y_5)} (-1)) \simeq \tor_5(\Bbbk,R/J).
\] 
To complete the argument, we will show that $\beta_{5,6}(R/(J :y_5)) \neq 0$.
Since $(J:y_5) = (x_1, y_2, x_3, x_4x_5, x_5x_6)$, this is equivalent to
$\beta_{2,3}(R/(x_4x_5,x_5x_6)) \neq 0$, which is true.

We showed so far that $c \leq 6$. Now, if $c < 5$, $\reg R/I < 3$. Hence
$c=5$ or $c=6$. As we noted in Discussion~\ref{disc:FirstSkipInml} that
$G$, therefore, is the suspension of the path or the cycle in five or six
vertices or of the $6$-cycle.

\eqref{qpureCMbipAtMostSixVert} $\implies$ \eqref{qpureCMbipRegGeqThree}:
If $G$ is the suspension of the path or the cycle on $c$ vertices, then
$\mathfrak d_G$ is such that every vertex is a source vertex or a sink
vertex, but not both, and that every source (respectively, sink) vertex in
$\mathfrak d_G$ is covered by (respectively, covers) at most two sink
(respectively, source) vertices. Hence $\reg R/I = \lceil \frac{c}2
\rceil$. Since $c=5$ or $c=6$ in our case, $\reg R/I = 3$. With this, 
$R/I$ has a quasi-pure resolution if and only if $m_4(I) = 6$, which we now
show. If on the other hand, $m_4(I) = 5$, then there exists $\sigma
\subseteq V$ and a partition $\sigma = \sigma_1 \bigsqcup \sigma_2$ (into
two sets, since $G$ is bipartite) such that $|\sigma| = 5$ and $G|_\sigma$
is a complete bipartite graph (Proposition~\ref{thm:ml-is-deg}). Recall
that $V = V_1 \bigsqcup V_2$ is the partition of the vertex set $V$ of $G$.
We may assume that $\sigma_i \subseteq V_i, i=1,2$. If $|\sigma_i| = 1$ for
any $i$, then $|\sigma| \leq 4$, because $\deg_G x \leq 3$ for all $x \in
V$. On the other hand, if, say, $|\sigma_1| \geq 2$, then $|\sigma_2| = 1$,
because otherwise, we would get a $4$-cycle in $G$, contradicting the fact
that $G$ has only a $6$-cycle, if any. Now, again, $|\sigma_1| \leq 3$, so
$|\sigma| < 5$. Hence $R/I$ has a quasi-pure resolution.
\end{proof}

We add, in passing, that the edge ideals $I$ of the suspension of
paths and cycles on four or fewer vertices have quasi-pure resolutions, but
this follows easily from the fact that $\reg R/I \leq 2$.

\section*{Acknowledgments}
The author thanks C. Huneke and J. Martin for helpful discussions.


\def\cfudot#1{\ifmmode\setbox7\hbox{$\accent"5E#1$}\else
  \setbox7\hbox{\accent"5E#1}\penalty 10000\relax\fi\raise 1\ht7
  \hbox{\raise.1ex\hbox to 1\wd7{\hss.\hss}}\penalty 10000 \hskip-1\wd7\penalty
  10000\box7}
\providecommand{\bysame}{\leavevmode\hbox to3em{\hrulefill}\thinspace}
\providecommand{\MR}{\relax\ifhmode\unskip\space\fi MR }
\providecommand{\MRhref}[2]{%
  \href{http://www.ams.org/mathscinet-getitem?mr=#1}{#2}
}
\providecommand{\href}[2]{#2}

\end{document}